\documentstyle[amssymb,amsfonts]{amsart}

\def\R{{\hbox{\bf R}}}

\def\I{{\hbox{\bf I}}}
\def\J{{\hbox{\bf J}}}

\font \roman = cmr10 at 10 true pt

\def\be#1{\begin{equation}\label{#1}}
\def\bas{\begin{align*}}
\def\eas{\end{align*}}
\def\bi{\begin{itemize}}
\def\ei{\end{itemize}}

\def\dist{{\hbox{\roman dist}}}

\newcommand{\qtil}{\tilde{q}}

\def\Z{{\hbox{\bf Z}}}
\def\eps{\varepsilon}
\newenvironment{proof}{\noindent {\bf Proof} }{\endprf\par}
\def \endprf{\hfill  {\vrule height6pt width6pt depth0pt}\medskip}
\def\emph#1{{\it #1}}
\def\textbf#1{{\bf #1}}

\theoremstyle{plain}
  \newtheorem{theorem}[subsection]{Theorem}

  \newtheorem{proposition}[subsection]{Proposition}
  
  \newtheorem{lemma}[subsection]{Lemma}
  \newtheorem{corollary}[subsection]{Corollary}

\theoremstyle{remark}

\theoremstyle{definition}
  \newtheorem{definition}[subsection]{Definition}

\include{psfig}

\begin{document}

\title{Endpoint multiplier theorems of Marcinkiewicz type}

\author{Terence Tao}
\address{Department of Mathematics, UCLA, Los Angeles CA 90095-1555}
\email{tao@@math.ucla.edu}

\author{Jim Wright}
\address{School of Mathematics, UNSW, Sydney Australia 2052}
\email{jimw@@maths.unsw.edu.au}


\begin{abstract}  We establish sharp $(H^1, L^{1,q})$ and local $(L \log^r L, L^{1,q})$ mapping properties for rough one-dimensional multipliers.  In particular, we show that the multipliers in the Marcinkiewicz multiplier theorem map $H^1$ to $L^{1,\infty}$ and $L \log^{1/2} L$ to $L^{1,\infty}$, and that these estimates are sharp.
\end{abstract}

\maketitle

\section{Introduction}

Let $m$ be a bounded function on $\R$, and let $T_m$ be the associated multiplier
$$ \widehat{T_m f}(\xi) = m(\xi) \hat f(\xi).$$
There are many multiplier theorems which give conditions under which $T_m$ is an $L^p$ multiplier. We will be interested in the mapping behaviour of $T_m$ near $L^1$.  Specifically, we address the following questions:

\begin{itemize}
\item For which $1 \leq q \leq \infty$ does $T_m$ map the Hardy space $H^1$ to the Lorentz space $L^{1,q}$?

\item We say that $T_m$ locally maps the Orlicz space $L \log^r L$ to $L^{1,q}$ if
$$ \| T_m f \|_{L^{1,q}(K)} \leq C_K \| f \|_{L \log^r L(K)}$$
for all compact sets $K$ and all functions $f$ on $K$.  For which $r \geq 0$ and $1 \leq q \leq \infty$ does $T_m$ locally map $L \log^r L$ to $L^{1,q}$?
\end{itemize}

Standard interpolation theory (see e.g. \cite{bergh:interp}) shows that if $T_m$ locally maps $L \log^r L$ to $L^{1,q}$, then it locally maps $L \log^{\tilde r}$ to $L^{1,\qtil}$ whenever $\qtil \leq q$ and $\tilde r \geq r + \frac{1}{\qtil} - \frac{1}{q}$.  Also, extrapolation theory (\cite{yano}, \cite{tao:factor}) shows that $T_m$ maps $L \log^r L$ to $L^1$ if and only if the $L^p$ operator norm of $T_m$ grows like $O((p-1)^{-r-1})$ as $p \to 1$.

Here and in the sequel, $\eta$ is an even bump function adapted to $\pm [1/2,4]$ which equals 1 on $\pm [3/4,3]$.

\begin{definition}  If $m$ is a symbol and $j$ is an integer, we define the \emph{$j^{th}$ frequency component} $m_j$ of $m$ to be the function
$$ m_j(\xi) = \eta(\xi) m(2^j \xi).$$
\end{definition}

We say that $T_m$ is a \emph{H\"ormander multiplier} if the frequency components $m_j$ are in the Sobolev space $L^2_{1/2+}$ uniformly in $j$.  These multipliers are Calder\'on-Zygmund operators and hence map $H^1$ to $L^1$ (and even to $H^1$), and $L^1$ to $L^{1,\infty}$; see e.g. \cite{stein:small}.  By interpolation one then sees
that $T_m$ locally maps $L \log^r L$ to $L^{1,q}$ whenever $r \leq 1/q$.

We now consider multipliers not covered by the H\"ormander theory.   We say that $T_m$ is a \emph{Marcinkiewicz multiplier} if the frequency components $m_j$ have bounded variation uniformly in $j$.  The Marcinkiewicz multiplier theorem (see e.g. \cite{stein:small}) shows that
$T_m$ is bounded on $L^p$.

Our first result characterizes the endpoint behaviour of Marcinkiewicz multipliers:

\begin{theorem}\label{marcin-thm} Marcinkiewicz multipliers map $H^1$ to $L^{1,\infty}$, and locally map $L \log^r L$ to $L^{1,q}$ whenever $r \geq \frac{1}{2} + \frac{1}{q}$.  Conversely, there exist Marcinkiewicz multipliers which do not map $H^1$ to $L^{1,q}$ for any $q < \infty$, and do not locally map $L \log^r L$ to $L^{1,q}$ for any $r < \frac{1}{2} + \frac{1}{q}$.
\end{theorem}

We can generalize the notion of a Marcinkiewicz multiplier as follows.  

\begin{definition}\cite{crdefs} Let $X$ denote the set of all functions of the form
$$ m  = \sum_I c_I \chi_I$$
where $I$ ranges over a collection of disjoint intervals in $\pm [1/2,4]$, and the $c_I$ are square summable coefficients:
\be{c-est}
(\sum_I |c_I|^2)^{1/2} \leq 1.
\end{equation}
Let $\overline X$ denote the Banach space generated by using the elements of $X$ as atoms; note that this space includes all functions of bounded variation on $\pm[1/2,4]$.  We say that $T_m$ is a \emph{$R_2$ multiplier} if the frequency components $m_j$ are in $\overline X$ uniformly in $j$.
\end{definition}

This class is more general than the Marcinkiewicz and H\"ormander classes.  In \cite{crdefs} it was established that $R_2$ multipliers are bounded on all $L^p$, $1 < p < \infty$.

We can extend the positive results of Theorem \ref{marcin-thm} as follows.

\begin{theorem}\label{rough-thm}  All the statements in Theorem \ref{marcin-thm} continue to hold for $R_2$ multipliers.
\end{theorem}

One can also show the $L^p$ norms of these multipliers grow like $\max(p,p')^{3/2}$ by converse extrapolation theorems (see \cite{tao:factor}).  This is sharp.  Theorem \ref{rough-thm} also has an easy corollary to multipliers of bounded $s$-variation as studied in \cite{crdefs}; we detail this in Section \ref{remarks}.

We now consider another multiplier class which is slightly smoother than the $R_2$ multiplier class.  

\begin{definition}\cite{seeger:compact}
Let $X'$ denote the set of all functions of the form
$$ m  = \sum_I c_I \psi_I$$
where $I$, $c_I$ are as in the definition of $X$, and the $\psi_I$ are $C^10$ bump functions adapted to $I$.  Let $\overline{X'}$ be the atomic Banach space generated by $X'$.  We say that $m$ is in $R_{1/2,2}^2$ if 
\be{whoop} \| \psi m(2^{-j} \cdot)  \|_{\overline{X'}} \lesssim 1
\end{equation}
for all integers $j$, where $\psi$ is a bump function adapted to $\pm [1/2,4]$ which equals 1 on $\pm [1,2]$.  We say that $T_m$ is a \emph{$R_{1/2,2}^2$ multiplier} if the frequency components $m_j$ are in $\overline{X'}$ uniformly in $j$.  
\end{definition}

This class was first studied in \cite{seeger:compact}; it contains the H\"ormander class, is contained in the $R_2$ class, and is not comparable with the Marcinkiewicz class.  In \cite{seeger:compact}, Theorem 2.2 the $R_{1/2,2}^2$ multipliers were shown to map $H^1$ to $L^{1,\infty}$; we can improve this to

\begin{theorem}\label{smooth-thm}  $R_{1/2,2}^2$ multipliers map $H^1$ to $L^{1,2}$, and locally map $L \log^r L$ to $L^{1,q}$ whenever $r \geq \max(\frac{1}{2}, \frac{1}{q})$.  Conversely, there exist $R_{1/2,2}^2$ multipliers which do not map $H^1$ to $L^{1,q}$ for any $q < 2$, and do not map $L \log^r L$ to $L^{1,q}$ whenever $r < \max(\frac{1}{2}, \frac{1}{q})$.
\end{theorem}

The converse extrapolation theorem in \cite{tao:factor} thus shows that these operators have an $L^p$ operator norm of $O(\max(p,p'))$, and this is sharp.

Thus, to summarize our main results, $R_2$ multipliers map both $H^1$ and $L \log L^{1/2}$ to $L^{1,\infty}$, while the smoother $R_{1/2,2}^2$ multipliers map both $H^1$ and $L \log L^{1/2}$ to $L^{1,2}$, with all exponents being best possible.

From the classical study \cite{hirschman:multiplier} of the multipliers
\be{mab}
m(\xi) = \frac{e^{i|\xi|^\alpha}}{(1+|\xi|^2)^{\beta/2}}
\end{equation}
it is known that the condition \eqref{c-est} cannot be replaced with a weaker $l^q$ condition, $q > 2$, if the intervals $I$ are the same size.  However, even if the intervals are different sizes, one still cannot relax this condition, as the following result shows.

\begin{definition}\cite{seeger:compact}  For any $1 \leq q \leq \infty$, let $X'_q$ be defined as in $X'$ but with \eqref{c-est} replaced by
$$ (\sum_k (\sum_{I: |I| \sim 2^k} c_I^2)^{q/2})^{1/q} \leq 1.$$
Let $\overline{X'_q}$ be the atomic Banach space generated by $X'_q$.  We say that $T_m$ is a \emph{$R_{1/2,q}^2$ multiplier} if the frequency components $m_j$ are in $\overline{X'_q}$ uniformly in $j$.
\end{definition}

\begin{theorem}\label{exotic}  For any $q > 2$, there exist $R_{1/2,q}^2$ multipliers which are unbounded on $L^p$ for $|\frac{1}{2} - \frac{1}{p}| > \frac{1}{q}$.  In particular, there are no mapping properties near $L^1$.
\end{theorem}

One can obtain positive $(L^p,L^p)$ or $(L^p,L^{p,2})$ mapping results when $2 < q \leq \infty$ for these operators by complex interpolation between Theorem \ref{smooth-thm} and trivial $L^2$ estimates (cf. \cite{feff:thesis}), but we shall not do so here. 

The space $H^1$ has of course appeared countless times in endpoint multiplier theory, but the appearance of the Orlicz space $L \log^{1/2} L$ space is more unusual.  This space first appeared in work of Zygmund \cite{zygmund}, who showed the inequality
\be{lacun} (\sum_{j = 0}^\infty |\hat f(2^j)|^2)^{1/2} \lesssim \|f\|_{L \log^{1/2} L}
\end{equation}
for all $f$ on the unit circle $S^1$.  This inequality can be viewed as a rudimentary prototype of the multiplier theorems described above (indeed, one can derive \eqref{lacun} from either of the above theorems by transplanting the results to the circle, and considering multipliers supported on the dyadic frequencies $2^j$).  As we shall see in Section \ref{sqfn}, the space $L \log^{1/2} L$ is in fact very similar to the Hardy space $H^1$ in that it has an associated square function which is integrable.

The space $L^{1,2}$ has appeared in recent work of Seeger and Tao \cite{ts:l12}.  Very roughly speaking, just as the space $L^{1,\infty}$ is natural for maximal functions and $L^1$ is natural for sums, the space $L^{1,2}$ is natural for certain square functions.  A concrete version of this principle appears in Lemma \ref{l12}.

This paper is organized as follows.  After some notational preliminaries we detail the negative results to the above Theorems in Section \ref{negative}.  
In Section \ref{sqfn} and the Appendix we show how both $H^1$ and $L \log^{1/2} L$ functions are associated with an integrable square function.  In Sections \ref{positive}, \ref{li-sec}, \ref{l2-sec} we then show how control of this square function leads to $L^{1,2}$ and $L^{1,\infty}$ multiplier estimates.  Finally, we discuss the $V_q$ class in Section \ref{remarks}.

This work was conducted at UNSW.  The authors thank Gerd Mockenhaupt and Andreas Seeger for useful comments.  The first author is supported by NSF grant DMS-9706764.

\section{Notation}

We use $C$ to denote various constants, and $A \lesssim B$, $A = O(B)$, or ``$B$ majorizes $A$'' to denote the estimate $A \leq CB$.  We use $A \sim B$ to denote the estimate $A \lesssim B \lesssim A$.

Here and in the sequel, $\Delta_j$ denotes the Littlewood-Paley multiplier with symbol $\eta(2^{-j} \cdot)$, where $\eta$ is as in the introduction.  For integers $j$, we use $\phi_j$ to denote the weight function
\be{phij-def} \phi_j(x) = 2^j (1 + 2^{2j} |x|^2)^{-3/4}.
\end{equation}
Similarly, for intervals $I$ we use $\phi_I$ to denote the weight
\be{phii-def} \phi_I(x) = |I| (1 + |I|^2 |x|^2)^{-3/4}.
\end{equation}
These weights are thus smooth and decay like $|x|^{-3/2}$ at infinity.  Many quantities in our argument will be controlled using the $\phi_j$, $\phi_I$; the reason why the decay is so weak is because we are forced at one point to use the Haar wavelet system, which has very poor moment conditions.  (The exact choice of $3/2$ has no significance, any exponent strictly between 1 and 2 would have sufficed).

\section{Negative results}\label{negative}

In this section we detail the counter-examples which yield the negative results stated in the introduction. In all of these examples $N$ is a large integer which will eventually be sent to infinity, $(e_j)_{j \in \Z}$ is the standard basis of $l^2(\Z)$, and $\psi$ is a non-negative even bump function supported on $\{ |\xi| \ll 1 \}$ which equals 1 at the origin and has a non-negative Fourier transform.  Some of our counter-examples will be vector-valued, but one can obtain scalar-valued substitutes by replacing $e_j$ with randomized signs $\eps_j = \pm 1$ and using the Lorentz-space version of Khinchin's inequality; we omit the details.

\subsection{Marcinkiewicz multipliers and $R_2$ multipliers need not map $H^1$ to $L^{1,q}$ for any $q < \infty$.}

Consider the symbol
\be{m0-def}
m_0(\xi) = \chi_{[1,\infty)}(\xi)\psi(\xi - 1).
\end{equation}
The convolution kernel $\hat{m_0}$ of this function is bounded for $|x| \lesssim 1$, and can be estimated via stationary phase as
\be{scoot} \hat{m_0}(x) = e^{2\pi i x}/x + O(|x|^{-2})
\end{equation}
for $|x| \gg 1$.  If we then test this multiplier against a bump function $f$ with $\hat f(0) = 0$ and $\hat f(1) \neq 0$, we see that $f$ is in $H^1$, but $|T_{m_0} f(x)| \sim 1/x$ as $|x| \to \infty$, so $T_{m_0} f$ is not in $L^{1,q}$ for any $q < \infty$.

\subsection{Marcinkiewicz multipliers and $R_2$ multipliers need not locally map $L \log^r L$ to $L^{1,q}$ for any $r < \frac{1}{2} + \frac{1}{q}$.}

Define the vector-valued multiplier 
$$ m_N(\xi) = \sum_{j = 0}^N e_j m_0(\xi/2^j)$$
where $m_0$ is defined in \eqref{m0-def}; this multiplier satisfies the requirements of both Theorems.  

By testing $T_{m_N}$ against a function $f$ whose Fourier transform is a bump function which equals 1 on $[-2^N, 2^N]$ and is adapted to a slight dilate of this interval, (so that $\|f\|_{L \log^r L} \sim N^{1/r}$) we see that we must have
$$ \| \hat{m_N} \|_{L^{1,q}([0,1])} \lesssim N^{1/r}$$
in order for $T_{m_N}$ to locally\footnote{Strictly speaking, $f$ is not quite compactly supported, but the error incurred because of this is extremely rapidly decreasing in $N$ and can be easily dealt with.} map $L \log^r L$ to $L^{1,q}$.
However, by \eqref{scoot} we have
$$ |\hat{m_N}(x)| \sim \frac{\log(1/|x|)^{1/2}}{|x|}$$
for $2^N \ll |x| \ll 1$, and the necessary condition $r < \frac{1}{2} + \frac{1}{q}$ follows by a routine computation. 

\subsection{$R^{2,1/2}_2$ multipliers need not map $H^1$ to $L^{1,q}$ for any $q < 2$.}

We use the multiplier
$$ m'_N(\xi) = N^{-1/2} \sum_{j=0}^N \psi(2^j(\xi-1)-1).$$
This multiplier is in the class of Theorem \ref{smooth-thm}.  Now suppose for contradiction that $T_{m'_N}$ mapped $H^1$ to $L^{1,q}$.  Since $m'_N$ is supported in a single dyadic scale, we may factor $T_{m'_N} = T_{m'_N} S_0$ where $S_0$ is a Littlewood-Paley projection to frequencies $|\xi| \sim 1$.  From the Littlewood-Paley square-function characterization we see that $S_0$ maps $H^1$ to $L^1$, hence $T_{m'_N}$ maps $L^1$ to $L^{1,q}$.  In particular, the kernel $\widehat{m'_N}$ must be in $L^{1,q}$.  However, a computation shows that
$$ |\widehat{m'_N}(x)| \lesssim \frac{N^{-1/2}}{|x|}$$
for $1 \ll |x| \ll 2^N$, which contradicts the assumption that $q < 2$.

\subsection{$R_{1/2,2}^2$ multipliers need not locally map $L \log^r L$ to $L^{1,q}$ for any $r < \frac{1}{2}$.}

We consider the vector-valued multiplier
$$ m''_N(\xi) = \sum_{j=0}^N e_j \psi(\xi - 2^j);$$
this is a multiplier in the class of Theorem \ref{smooth-thm}.  By repeating the argument with the $m_N$ multipliers, we must have
$$ \| \hat{m''_N} \|_{L^{1,q}([0,1])} \lesssim N^{1/r}.$$
However, a computation shows that
$$ |\hat{m''_N}(x)| \sim \sqrt{N}$$
for $|x| \ll 1$, and this contradicts the assumption $r < \frac{1}{2}$.  

\subsection{$R_{1/2,2}^2$ multipliers need not locally map $L \log^r L$ to $L^{1,q}$ for any $r < \frac{1}{q}$.}

We consider the Hilbert transform $H$, which of course is of the class in Theorem \ref{smooth-thm}, and test it against the function 
$f = 2^N \chi_{[0,2^{-N}]}$.  Clearly $f$ has a $L \log^r L$ norm of $N^r$ but the Hilbert transform of this function has a local $L^{1,q}$ norm of about $N^{1/q}$, hence the claim.

\subsection{$R_{1/2,q}^2$ multipliers need not be bounded on $L^p$ for $|\frac{1}{2} - \frac{1}{p}| > \frac{1}{q}$.}

By duality it suffices to show unboundedness when $\frac{1}{p} - \frac{1}{2} > \frac{1}{q}$.

We define the vector-valued multiplier
$$ m'''_N(\xi) = N^{-1/q} \sum_{j=N/100}^{N/10} e_j \psi(2^j(\xi - \frac{j}{N}));$$
this multiplier is in the class of Theorem \ref{exotic}.  We test this against the function
$$ f(x) = \sum_{|k| < 2^N} \psi(x - Nk).$$
We expand
$$ T_{m'''_N} f(x) = N^{1/q} \sum_{j=N/100}^{N/10} e_j \sum_{|k| < 2^N}
\int \psi(x-y-Nk) e^{2\pi i j y/N} 2^{-j} \hat \psi(2^{-j}y)\ dy.$$
Making the change of variables $y \to y - Nk$, this becomes
$$ N^{1/q} \sum_{j=N/100}^{N/10} e_j \sum_{|k| < 2^N}
\int \psi(x-y) e^{2\pi i j y/N} 2^{-j} \hat \psi(2^{-j}(y+Nk))\ dy.$$
The function $e^{2\pi i j y/N}$ has real part bounded away from zero, so
$$ |T_{m'''_N} f(x)| \sim
N^{-1/q} (\sum_{j=N/100}^{N/10} 
(\int \psi(x-y) 2^{-j} \sum_{|k| \leq K} \hat \psi(2^{-j}(y+Nk))\ dy)^2)^{1/2}.$$
If $|x| \ll 2^N$, then $|y| \ll 2^N$ and the inner sum is $\sim 2^j/N$ (note that $N2^N \gg 2^j \gg N$).  Thus we have
$$ |T_{m'''_N} f(x)| \sim
N^{1/q} (\sum_{j=N/100}^{N/10} 
(\int N^{-1} \psi(x-y)\ dy)^2)^{1/2} \sim N^{-1/q - 1/2}$$
for $|x| \ll 2^N$.  Thus
$$ \| T_{m'''_N} f\|_p \gtrsim N^{-1/q - 1/2} 2^{N/p}.$$
On the other hand, an easy computation shows
$$ \|f\|_p \sim N^{-1/p} 2^{N/p},$$
which demonstrates unboundedness when $\frac{1}{p} - \frac{1}{2} > \frac{1}{q}$.

\section{The spaces $H^1$ and $L \log^{1/2} L$.}\label{sqfn}

Our positive results involve the spaces $H^1$ and $L \log^{1/2} L$.  As is well known, $L \log^{1/2} L$ functions are in general not in $H^1$ and thus do not have an integrable Littlewood-Paley square function.  However, there is a substitute square function for these functions which are indeed integrable, which is why all our results for $H^1$ also extend to $L \log^{1/2} L$.  More precisely:

\begin{proposition}\label{cont}  Let $f$ be a function which is either in the unit ball $H^1(\R)$, or in the unit ball of $L \log^{1/2} L([-C,C])$ and with mean zero.  Then there exists non-negative functions $F_j$ for each integer $j$ such that we have the pointwise estimate
\be{fj-support}
|\Delta_j f(x)| \lesssim F_j * \phi_j(x)
\end{equation}
for all $j \in \Z$ and $x \in \R$, and the square function estimate
\be{fj-norm}
\| (\sum_j |F_j|^2)^{1/2} \|_1 \lesssim 1
\end{equation}
\end{proposition}

This proposition is easy to prove when $f$ is in $H^1$.  Indeed, one simply chooses $F_j = |\tilde \Delta_j f|$, where $\tilde \Delta_j$ is a slight enlargement of $\Delta_j$ such that $\Delta_j = \Delta_j \tilde \Delta_j$.  
The claim \eqref{fj-support} follows from pointwise control on the kernel of $\Delta_j$, while \eqref{fj-norm} follows from the square function characterization of $H^1$.

The corresponding claim for $L \log^{1/2} L$ is much more delicate.  We remark that this claim implies Zygmund's inequality \eqref{lacun}.  To see this, we first observe that we may assume $f$ satisfies the conditions of the above Proposition, in which case $\hat f(2^j)$ can  be estimated by $\| \Delta_j f\|_1 \lesssim \|F_j\|_1$.  The claim then follows from \eqref{fj-norm} and the Minkowski inequality
$$ (\sum_j \|F_j\|_1^2)^{1/2} \leq \| (\sum_j |F_j|^2)^{1/2} \|_1.$$
The same argument shows that $L \log L^{1/2}$ cannot be replaced by any weaker Orlicz norm.  However, the Proposition is substantially stronger than Zygmund's inequality.

As an example of the Proposition, let $f = 2^N N^{-1/2} \psi_N$, where $N$ is a large integer and $\psi_N$ is a bump function of mean zero adapted to the interval $[-2^{-N},2^{-N}]$.  This function is normalized in $L \log^{1/2} L$ and has mean zero, but is not in $L^1$.  Indeed, if one lets $F_j = |\tilde \Delta_j f|$ as before, then for each $1 \ll j \ll N$, $F_j$ is comparable to $2^j N^{-1/2} \psi_j$ on the interval $[-2^{-j}, 2^{-j}]$, and is rapidly decreasing outside of this interval.  From this we see that the left hand side of \eqref{fj-norm} is too large (about $N^{1/2}$).  The problem here is that the functions $F_j$ have very different supports, and so their contributions to \eqref{fj-norm} add up in $l^1$ rather than $l^2$.  To get around this we can redistribute the mass of the $F_j$, setting $F_j = 2^N N^{-1/2} \chi_{[-2^{-N}, 2^{-N}]}$ for each $1 \ll j \ll N$; one verifies that \eqref{fj-support} is still satisfied, and that \eqref{fj-norm} is now satisfied because the $F_j$ are summing in $l^2$ rather than $l^1$.  (The frequencies $j \leq 1$ or $j \geq N$
can be handled by the original assignment $F_j = |\Delta_j f|$ without difficulty).

To handle the general case we shall follow a similar philosophy, namely that each $F_j$ shall be a redistribution of $|\Delta_j f|$, whose supports overlap so much that their contributions to \eqref{fj-norm} are summed in $l^2$ rather than $l^1$.  To do this for general functions $f$ we will use a delicate recursive algorithm.  In order to control the error terms in this algorithm we shall be forced to move to the dyadic (Haar wavelet) setting, and also to reduce $f$ to a characteristic function.

The argument is somewhat lengthy, and the methods used are not needed anywhere else in the paper.  Because of this, we defer the argument to an Appendix, and proceed to the key estimate in the proofs of Theorems \ref{rough-thm}, \ref{smooth-thm} in the next section.

\section{Positive results: the main estimate}\label{positive}

In this section we summarize the main estimate we will need to prove in order to achieve the positive results in Theorems \ref{rough-thm} and \ref{smooth-thm}.  (The positive results in Theorem \ref{marcin-thm} follow immediately from those in Theorem \ref{rough-thm}).

By interpolation with the trivial $L^2$ boundedness results coming from Plancherel's theorem, it suffices to show that the operators in Theorem \ref{rough-thm} map $H^1$ and $L \log^{1/2} L$ to $L^{1,\infty}$, and the operators in Theorem \ref{smooth-thm} map $H^1$ and $L \log^{1/2} L$ to $L^{1,2}$.

We will use two key results to obtain these boundedness properties.  The first is the square function estimate obtained above in Proposition \ref{cont}.  The second is an endpoint multiplier result associated to an arbitrary collection of intervals, which we now state.

\begin{proposition}\label{main}  Let $N \geq 1$ be an integer, and let $\{I\}$ be a collection of intervals in $\R$ which overlap at most $N$ times in the sense that
\be{overlap} \| \sum_{I} \chi_I \|_\infty \leq N.
\end{equation}
For each $I$, we assign a function $f_I$, a non-negative function $F_I$, and a multiplier $T_{m_I}$ with the following properties.
\bi
\item For each $I$, $m_I$ is supported on $I$, there exists a $\xi_I \in I$ such that the symbol $m_I(\cdot + \xi_I)$ is a standard symbol of order 0 in the sense of e.g. \cite{stein:large}.
\item For any $I \in \I$ and $x \in \R$ we have the pointwise estimate
\be{point}
|f_I(x)| \lesssim F_I(x) * \phi_I(x)
\end{equation}
where $\phi_I$ was defined in \eqref{phii-def}.
\ei
Then we have
\be{linfty-est}
\| \sum_{I } T_{m_I} f_I \|_{L^{1,\infty}} \lesssim N^{1/2} \| (\sum_{I} |F_I|^2)^{1/2} \|_1.
\end{equation}
If we strengthen the condition on $m_I$ and assume that the $m_I$ are actually bump functions adapted to $I$ uniformly in $I$, then we may strengthen \eqref{linfty-est} to
\be{l2-est}
\| \sum_{I} T_{m_I} f_I \|_{L^{1,2}} \lesssim N^{1/2} \| (\sum_{I} |F_I|^2)^{1/2} \|_1.
\end{equation}
\end{proposition}

We will prove this proposition in Sections \ref{li-sec}, \ref{l2-sec}.  For now, we see how this proposition and Proposition \ref{cont} imply the desired mapping properties on $R_2$ and $R^2_{1/2,2}$ multipliers.

Let us first make the preliminary reduction that to prove the $L \log^{1/2} L$ local mapping properties on $T_m$ it suffices to prove global estimates on $T_m f$ assuming that $f$ is supported in $[0,1]$, is normalized in $L \log^{1/2} L$, and has mean zero.  The normalization to $[0,1]$ follows from dilation and translation invariance; the mean zero assumption comes by subtracting off a bump function and observing from the $L^2$ theory that $T_m$ applied to a bump function is locally in $L^2$, hence locally in $L^{1,\infty}$ and $L^{1,2}$.

Our task is now to show that any $f$ satisfying either of the conditions in Proposition \ref{cont}, we have
\be{mi} \|T_m f\|_{L^{1,\infty}} \lesssim 1
\end{equation}
for $R_2$ multipliers and
\be{m2}
 \|T_m f\|_{L^{1,2}} \lesssim 1
\end{equation}
for $R^2_{1/2,2}$ multipliers.

Fix $f$, and let $F_j$ be as in Proposition \ref{cont}.  We first prove \eqref{mi}. We may assume without loss of generality that $m$ is supported in $\bigcup_{j \hbox{ even}} [2^j,2^{j+1}]$ (The case of odd $j$ is similar and is omitted).
By a limiting argument we may assume that only finitely many of the frequency components $m_j$ are non-zero for even $j$.  By a further limiting argument we may assume that each $m_j$ for even $j$ is a rational linear combination of elements in $X$, e.g. $m_j = \sum_{i=1}^{N_j} \alpha_{j,i} m_{j,i}$ where the $m_{j,i}$ are uniformly in $X$ and the $\alpha_{j,i}$ are non-negative rational numbers.  By placing the rational $\alpha_{j,i}$ under a common denominator $N$, and repeating each $m_{j,i}$ with a multiplicity equal to $N \alpha_{j,i}$, we may thus write
$$ m = \frac{1}{N} \sum_{i=1}^N m^{(i)}$$
where the frequency components $m^{(i)}_j$ are uniformly in $X$ for even $j$.  In particular, this implies that
$$ m = \sum_I c_I \chi_I$$
where each interval $I$ belongs to $[2^{j_I}, 2^{j_I+1}]$ for some even $j_I$, the intervals $I$ satisfy \eqref{overlap}, and
\be{c-bound} \sum_{I: j_I = j} c_I^2 \lesssim N^{-1}
\end{equation}
for each $j$.  We may assume that $|I| \ll 2^{j_I}$ for all $I$.
We split $\chi_I$ as
\be{chi-decomp} \chi_I(\xi) = \psi_I \psi^l_I H(\xi - \xi^l_I) + \psi_I \psi^r_I H(\xi^r_I - \xi)
\end{equation}
where $H = \chi_{(0,\infty)}$ is the Heaviside function, $\xi^l_I$ and $\xi^r_I$ are the left and right endpoints of $I$, and $\psi^l_I$, $\psi^r_I$, $\psi_I$ are bump functions adapted to $[\xi_l - |I|, \xi_l + |I|]$, $[\xi_r - |I|, \xi_r + |I|]$, and $5I$ respectively.  

We thus need to prove
$$
\| \sum_I c_I T_{\psi_I} T_{\psi^l_I H(\cdot - \xi^l_I)} f \|_{L^{1,\infty}} \lesssim 1,$$
together with the analogous estimate with the $l$ index replaced by $r$.  We show the displayed estimate only, as the other estimate is proven similarly.

Write $m_I = \psi^l_I H(\cdot - \xi^l_I)$, $\xi_I = \xi^l_I$, $f_I = c_I T_{\psi_I} f$, and $F_I = |c_I| F_{j_I}$.  The estimate \eqref{point} follows from eqref{fj-support}, the identity $T_{\psi_I} = T_{\psi_I} \Delta_{j_I}$ and kernel estimates on $T_{\psi_I}$.  Applying \eqref{linfty-est} we thus see that
$$
\| \sum_I c_I T_{\psi_I} T_{\psi^l_I H(\cdot - \xi^l_I)} f \|_{L^{1,\infty}} \lesssim N^{1/2} \| (\sum_{I} |F_I|^2)^{1/2} \|_1.
$$
The claim then follows from the definition of $F_I$, \eqref{c-bound}, and \eqref{fj-norm}.  This proves \eqref{mi}

The proof of \eqref{m2} is similar, but with $\chi_I$ replaced by a bump function $\tilde \psi_I$ adapted to $I$.  The only change is that the splitting \eqref{chi-decomp} is replaced by $\tilde \psi_I = \psi_I \tilde \psi_I$, where $\psi_I$ is a bump function adapted to $5I$ which equals 1 on $I$, and that \eqref{l2-est} is used instead of \eqref{linfty-est}.

It remains only to prove \eqref{linfty-est} and \eqref{l2-est}.  This shall be done in the next two sections.

\section{Proof of \eqref{linfty-est}}\label{li-sec}

Fix $I$, $N$, $f_I$, $F_I$, $m_I$; we may assume by limiting arguments that the collection of $I$ is finite.  From \eqref{point} we can find bounded functions $a_I$ for each $I \in \I$ such that
$$ f_I = a_I (F_I * \phi_I).$$
Our task is then to show that
$$
| \{ |\sum_{I} T_{m_I} (a_I (F_I * \phi_I))| \gtrsim \alpha \} |
\lesssim \alpha^{-1} N^{1/2} \| F \|_1$$
where $F$ denotes the vector $F = (F_I)_{I \in \{I\}}$.

We now perform a standard vector-valued Calder\'on-Zygmund decomposition
on $F$ at height $N^{-1/2} \alpha$ as
$$ F = g + \sum_J b_J$$
where $g = (g_I)_{I \in \I}$ satisfies the $L^2$ estimate
\be{g-bound}
\| g \|_2^2 \lesssim N^{-1/2} \alpha \| F \|_1,
\end{equation}
while the bad functions $b_J$ are supported on $J$, satisfy the moment condition $\int_J b_J = 0$, and the $L^1$ estimate
$$ \| b_J \|_1 \lesssim N^{-1/2} \alpha |J|.$$
Finally, the intervals $J$ satisfy
$$ \sum_J |J| \lesssim \alpha^{-1} N^{1/2} \| F \|_1.$$

Consider the contribution of the good function $g$.  By Chebyshev, it suffices to prove the $L^2$ estimate
\be{g-bit} \|\sum_{I} T_{m_I} (a_I (g_I * \phi_I)) \|_2^2 \lesssim
\alpha N^{1/2} \| (\sum_{I} |F_I|^2)^{1/2} \|_1.
\end{equation}
From Plancherel, the overlap condition on the $I$, and Cauchy-Schwarz, we have
the basic inequality
\be{basic}
\| \sum_{I} T_{m_I} h_I \|_2^2 \leq N \sum_I \|T_{m_I} h_I\|_2^2
\end{equation}
for any $h_I$.  We may thus estimate the left-hand side of \eqref{g-bit} by
\bas
N \sum_I \|T_{m_I} (a_I (g_I * \phi_I)) \|_2^2 &\lesssim
N \sum_I \|a_I (g_I * \phi_I) \|_2^2 \\
&\lesssim N \sum_I \|g_I * \phi_I \|_2^2 \\
&\lesssim N \sum_I \|g_I \|_2^2 \\
&\lesssim N N^{-1/2} \alpha \| (\sum_{I} |F_I|^2)^{1/2} \|_1
\end{align*}
as desired.

It remains to deal with the bad functions $b_J$.  It suffices to show that
$$
| \{ |\sum_{I} \sum_J T_{m_I} (a_I (b_{J,I} * \phi_I))| \gtrsim \alpha \} |
\lesssim \sum_J |J|.$$

From uncertainty principle heuristics we expect the contribution of the case $|I| |J| \leq 1$ to be easy.  Indeed, this case can be treated almost exactly like the good function $g$.  As before, it suffices to show the $L^2$ estimate
$$
\| \sum_{I,J: |I||J| \leq 1} T_{m_I} (a_I (b_{J,I} * \phi_I))\|_2^2 \lesssim
\alpha^2 \sum_J |J|.$$
By repeating the previous calculation, the left-hand side is majorized by
$$ N \sum_I \| \sum_{J: |I||J| \leq 1} b_{J,I} * \phi_I \|_2^2.$$
From the triangle inequality, it thus suffices to show that
$$ \sum_I \| \sum_{J: |I||J| = 2^{-m}} b_{J,I} * \phi_I \|_2^2
\leq 2^{-2m} N^{-1} \alpha^2 \sum_J |J|$$
for all $m \geq 0$.  This in turn follows if we can show
\be{unc} \sum_{I: |I| = 2^{-m-j}} \| \sum_{J: |J| = 2^j} b_{J,I} * \phi_I \|_2^2
\leq 2^{-2m} N^{-1} \alpha^2 \sum_{J: |J| = 2^j} |J|
\end{equation}
for all $m \geq 0$ and $j \in \Z$.

Fix $m$, $j$, and observe from \eqref{phij-def} that $\phi_I = \phi_{-m-j}$.  By moving the $I$ summation inside the norm, we can estimate the left-hand side of \eqref{unc} by 
$$\| \sum_{J: |J| = 2^j} b_J * \phi_{-m-j} \|_2^2$$
where $*$ is now a vector-valued convolution.  From the normalization and moment condition on $b_J$ we have
$$ b_J * \phi_{-m-j} \lesssim N^{-1/2} \alpha \chi_J * \phi_{-m-j}.$$
Inserting this into the previous, the claim then follows from Young's inequality and the $L^1$ normalization of the $\phi_{-m-j}$.

It remains to treat the case $|I||J| > 1$.  We split
$$ b_{J,I} * \phi_I = \chi_{2J} (b_{J,I} * \phi_I) + (1-\chi_{2J}) (b_{J,I} * \phi_I).$$
The contribution of the latter terms can be dealt with in a manner similar to that of the $|I| |J| \leq 1$ case.  As before, it suffices to show the $L^2$ estimate
$$
\| \sum_{I,J: |I||J| > 1} T_{m_I} (a_I (1-\chi_{2J}) (b_{J,I} * \phi_I))\|_2^2 \lesssim \alpha^2 \sum_J |J|.$$
As before, the left-hand side is majorized by
\be{temp} N \sum_I \| \sum_{J: |I||J| > 1} (1 - \chi_{2J}) (b_{J,I} * \phi_I) \|_2^2.
\end{equation}
A computation shows the pointwise estimate
$$ |(1 - \chi_{2J}) (b_{J,I} * \phi_I)| \lesssim \| b_{J,I}\|_1 |J|^{-1} (M \chi_J)^{3/2}.$$
(In fact there is an additional decay if $|I| |J|$ is large, but we shall not exploit this).  Inserting this estimate into \eqref{temp} and moving the $I$ summation back inside, we can majorize \eqref{temp} by
$$ N \| (\sum_I |\sum_J \|b_{J,I}\|_1 |J|^{-1} (M \chi_J)^{3/2}|^2)^{1/2} \|_2^2.$$
Using the triangle inequality for $l^2$ we may move the $I$ square-summation inside the $J$ summation.  If one then applies Minkowski's inequality
\be{mink}
(\sum_I \| b_{J,I}\|_1^2)^{1/2} \leq \| b_J \|_1 \lesssim N^{-1/2} \alpha |J|
\end{equation}
we can thus majorize \eqref{temp} by
$$ \alpha^2 \| \sum_J (M \chi_J)^{3/2}\|_2^2.$$
The claim then follows from the Fefferman-Stein vector-valued maximal inequality \cite{feff:max}.

It remains to show that
\be{new}
| \{ |\sum_{I,J: |I||J| > 1} T_{m_I} B_{J,I}| \gtrsim \alpha \} |
\lesssim \sum_J |J|
\end{equation}
where
$$ B_{J,I} = a_I \chi_{2J} (b_{J,I} * \phi_I).$$
For future reference we note from \eqref{mink} that the $B_{J,I}$ are supported on $2J$ and satisfy
\be{mink2}
\sum_I \| B_{J,I}\|_1^2 \lesssim N^{-1} \alpha^2 |J|^2
\end{equation}
for all $J$.

For each $I$, $J$ in \eqref{new}, let $P_{J,I}$ be a multiplier whose symbol is a bump function which equals $1$ on the interval $[\xi_I - |J|^{-1}, \xi_I + |J|^{-1}]$, and is adapted to a dilate of this interval.  We split 
$$ T_I = T_I P_{J,I} + Q_{J,I}$$
where $Q_{J,I} = T_I (1 - P_{J,I})$.  The point is that even though the kernel of $T_I$ decays very slowly, the operators $P_{J,I}$ and $Q_{J,I}$ have kernels which are essentially supported on an interval of width $|J|$.

We first consider the contribution of the $T_I P_{J,I}$.  It suffices as before to prove an $L^2$ estimate:
\be{pij}
\| \sum_{I,J: |I||J| > 1} T_{m_I} P_{J,I} B_{J,I}\|_2^2 \lesssim
\alpha^2 \sum_J |J|.
\end{equation}
By \eqref{basic} again, the left-hand side of \eqref{pij} is majorized by
$$ N \sum_I \| \sum_{J: |I| |J| > 1} P_{J,I} B_{J,I} \|_2^2.$$
From kernel estimates on $P_{I,J}$ we have the pointwise estimates
$$
|P_{J,I} B_{J,I}| \lesssim \| B_{J,I}\|_1 |J|^{-1} (M \chi_J)^{3/2}.$$
The contribution of the $T_I P_{J,I}$ is thus acceptable by repeating the arguments used to treat \eqref{temp}, and using \eqref{mink2} instead of \eqref{mink}.

It remains to consider the contribution of the $Q_{J,I}$.  For this final contribution we will not use $L^2$ estimates, but the more standard $L^1$ estimates outside an exceptional set:
$$ \| \sum_{I,J: |I||J| > 1} Q_{J,I} B_{J,I}\|_{L^1((\bigcup_J CJ)^c)} \lesssim \alpha \sum_J |J|.$$
By the triangle inequality it suffices to prove this for each $J$ separately:
$$ \| \sum_{I: |I||J| > 1} Q_{J,I} B_{J,I} \|_{L^1((CJ)^c)} \lesssim \alpha |J|.$$

By translation and scale invariance we may set $J = [0,1]$.  Let $\varphi$ denote a bump function which equals 1 on $[-1,1]$ and is adapted to $[-2,2]$.  Let $r_I$ denote the symbol
$$ r_I = q_{J,I} - q_{J,I} * \varphi,$$
where $q_{J,I}$ is the symbol of $Q_{J,I}$.  Observe that $Q_{J,I} B_{J,I} = T_{r_I} B_{J,I}$ outside of $CJ$.  Thus it suffices to show that
$$ \| \sum_{I: |I| > 1} T_{r_I} B_{J,I} \|_{L^1((CJ)^c)} \lesssim \alpha.$$
By H\"older's inequality it suffices to show the global weighted $L^2$ estimate
$$ \| x \sum_{I: |I| > 1} T_{r_I} B_{J,I}(x) \|_2 \lesssim \alpha.$$
By Plancherel, this becomes
$$ \| \sum_{I: |I| > 1} (r_I \widehat{B_{J,I}})' \|_2 \lesssim \alpha,$$
where the prime denotes differentiation.

The function $\widehat{B_{J,I}}$ is very smooth, in fact it satisfies the estimates
$$ \| \widehat{B_{J,I}} \|_{C^1} \lesssim \|B_{J,I}\|_1$$
for all $I$.  A computation using the construction of $Q_{J,I}$ and $r_I$ shows that the symbol $r_I$ satisfies the estimates
$$ |r_I(\xi)|, |r'_I(\xi)| \lesssim (1 + |\xi - \xi_I|)^{-10}$$
Combining these two estimates we see the pointwise estimate
$$ |(r_I \widehat{B_{J,I}})'| \lesssim \| B_{J,I} \|_1 (M \chi_{[\xi_I - 1, \xi_I + 1]})^2.$$
From the Fefferman-Stein vector-valued maximal inequality \cite{feff:max} it thus suffices to show that
$$ \| \sum_{I: |I| > 1} \| B_{J,I}\|_1 \chi_{[\xi_I-1,\xi_I+1]} \|_2 \lesssim \alpha.$$
However from \eqref{overlap} and the hypothesis $|I| > 1$ we see that the characteristic functions $\chi_{[\xi_I-1,\xi_I+1]}$ overlap at most $O(N)$ times at any given point.  The claim then follows from Cauchy-Schwarz and \eqref{mink2}.  This completes the proof of \eqref{linfty-est}.
\endprf

We remark that the one can modify this argument so that one does not need the full power of Proposition \ref{cont} in the $L \log^{1/2} L$ case, using a rescaled version of Zygmund's estimate \eqref{lacun} (for arbitrary lacunary frequencies, not just the powers of 2) as a substitute; we omit the details. 
On the other hand, the $(L \log^{1/2} L, L^{1,2})$ result in Proposition \ref{smooth-thm} seems to require the full strength of Proposition \ref{cont}.

\section{Proof of \eqref{l2-est}}\label{l2-sec}

We now prove \eqref{l2-est}.  As before we fix $I$, $N$, $m_I$, $f_I$, $F_I$, and assume that the collection of $I$ is finite.  We may also assume that the functions $F_I$ are smooth.

To prove \eqref{l2-est} it suffices to prove the stronger estimate
\be{l2-mod}
\| \sum_{I} T_{m_I} f_I \|_{L^{1,2}} \lesssim N^{1/2} \| (\sum_{I} |F_I * \phi_I|^2)^{1/2} \|_{L^{1,2}}.
\end{equation}

This is because of the following lemma, which illustrates the natural role of the Lorentz space $L^{1,2}$.

\begin{lemma}\label{l12} Let $I$ be an arbitrary collection of intervals, and $F_I$ an arbitrary collection of non-negative functions.  Then
$$ \|(\sum_I |F_I * \phi_I|^2)^{1/2} \|_{L^{1,2}} \lesssim
\| (\sum_I |F_I|^2)^{1/2} \|_1.$$
\end{lemma}

\begin{proof}  The desired estimate is the $p=2$ case of the more general estimate
$$ \|(\sum_I |F_I * \phi_I|^p)^{1/p} \|_{L^{1,p}} \lesssim
\| (\sum_I |F_I|^p)^{1/p} \|_1.$$
This estimate is trivial for $p=1$ by Young's inequality and the integrability of the $\phi_I$.  For $p=\infty$ the claim follows from the Hardy-Littlewood maximal inequality and the pointwise estimates
$$ |F_I * \phi_I(x)| \lesssim MF_I(x) \lesssim M (\sup_I F_I)(x).$$
The complex interpolation theorem of Sagher \cite{sagher} for Lorentz spaces then allows one to obtain the $p=2$ estimate.  Alternatively, one can interpolate manually by writing $F_I = |F| a_I$, where $|F| = (\sum_I |F_I|^2)^{1/2}$, and exploiting the Cauchy-Schwarz inequality
$$ |F_I * \phi_I(x)|^2 \leq  ((Fa_I^2) * \phi_j(x)) (|F| * \phi_j(x)) \lesssim |F|a_I^2 * \phi_I(x) M|F|(x)$$
and the H\"older inequality for Lorentz spaces \cite{oneil}
$$ \| (fg)^{1/2} \|_{L^{1,2}} \lesssim \|f\|_1^{1/2} \|g\|_{L^{1,\infty}}^{1/2}.$$
We omit the details.
\end{proof}

It remains to prove \eqref{l2-mod}.  Let $G$ denote the square function
$$ G = (\sum_{I} |F_I * \phi_I|^2)^{1/2};$$
note that $G$ is continuous from our a priori assumptions. It would be nice if the distributional estimate
$$
|\{ |\sum_{I} T_{m_I} f_I| \sim 2^j \}| \lesssim 
|\{ G \sim N^{-1/2} 2^j \}|$$
held for all $j$, as this easily implies \eqref{l2-mod}.  While this is not quite true, we are able to prove the substitute
\be{dist}
|\{ |\sum_{I} T_{m_I} f_I| \gtrsim 2^j \}| \lesssim 
2^{-2j} N \| \min(G, N^{-1/2} 2^j)\|_2^2
\end{equation}
for all $j$.  Indeed, if \eqref{dist} held, then we have
$$ 
2^j |\{ |\sum_{I} T_{m_I} f_I| \sim 2^j \}| \lesssim 
N^{1/2} \sum_s 2^{-|s|}
N^{-1/2} 2^{j+s} |\{ G \sim N^{-1/2} 2^{j+s} \}|.$$
the claim then follows by square-summing this in $j$, using the estimate
$$ \| F \|_{L^{1,2}} \sim (\sum_j (2^j |\{ F \sim 2^j\}|)^2)^{1/2}$$
and using Young's inequality.

It remains to prove \eqref{dist}.  Fix $j$, and consider the set $\Omega = \{ G > N^{-1/2} 2^j\}$. Since $G$ is continuous, $\Omega$ is an open set, and we may decompose it into intervals $\Omega = \bigcup_J J$ such that $G(x) = N^{-1/2} 2^j$ on the endpoints of $J$. Note that
\be{jit}
 \sum_J |J| = |\Omega| \leq 2^{-2j} N \| \min(G, N^{-1/2} 2^j) \|_2.
\end{equation}
We can therefore split
\be{spl}
\sum_I T_{m_I} f_I = \sum_I T_{m_I} (f_I\chi_{\Omega^c}) + \sum_{I,J: |I||J| \leq 1} T_{m_I} (f_I \chi_J) + \sum_{I,J: |I||J| > 1} T_{m_I} (f_I \chi_J).
\end{equation}

To treat the contribution of the first term in \eqref{spl} we use $L^2$ estimates.  By Chebyshev it suffices to show that
$$ \| \sum_{I} T_{m_I} (f_I\chi_{\Omega^c}) \|_2^2 \lesssim N \| \min(G, N^{-1/2} 2^j)\|_2^2.$$
However, by \eqref{basic} the left-hand side is majorized by
\bas
N \sum_I \| f_I \chi_{\Omega^c} \|_2^2 &=
N \| (\sum_I |f_I|^2)^{1/2} \chi_{\Omega^c} \|_2^2\\
&\lesssim N \| (\sum_I |F_I * \phi_I|^2)^{1/2} \chi_{\Omega^c} \|_2^2\\
&\leq N \| \min(G, N^{1/2} 2^j) \|_2^2
\end{align*}
as desired.

To treat the second term in \eqref{spl} we also use $L^2$ estimates.  As before, it suffices to show
\be{spl2} \| \sum_{I} T_{m_I} (\sum_{J: |I||J| \leq 1} f_I\chi_J) \|_2^2 \lesssim N \| \min(G, N^{-1/2} 2^j)\|_2^2.
\end{equation}
Using \eqref{basic} as before, we can majorize the left-hand side of \eqref{spl2} by
$$ 
N \sum_I \| \sum_{J: |I||J| \leq 1} (F_I * \phi_I)\chi_J \|_2^2.$$
Since the $J$ are all disjoint, we may re-arrange this as
$$ 
N \sum_J \sum_{I: |I||J| \leq 1} \|F_I * \phi_I \|_{L^2(J)}^2.$$
For each $J$ let $x_J^r$ be the right endpoint of $J$, so that $G(x_J^r) \leq N^{-1/2} 2^j$.  Now we exploit the assumption $|I||J| \leq 1$ to observe that
$$ |F_I * \phi_I(x)| \lesssim |F_I * \phi_I(x_J^r)|$$
for all $x \in J$.  Applying this to the previous, we can thus majorize \eqref{spl2} by
$$ 
N \sum_J |J| \sum_I |F_I * \phi_I(x_J^r)|^2 = N \sum_J |J| G(x_J^r)^2 \leq
2^{2j} \sum_J |J|.$$
The claim then follows from \eqref{jit}.

It remains to treat the third term in \eqref{spl}.  By Chebyshev and \eqref{jit} it suffices to prove an $L^1$ estimate outside the exceptional set $\bigcup_J CJ$:
$$ \| \sum_{I,J: |I||J| > 1} T_{m_I} (f_I\chi_J) \|_{L^1((\bigcup_J CJ)^c)} 
\lesssim 2^{j} \sum_J |J|.$$
By the triangle inequality it suffices to prove this for each $J$ separately:
$$ \| \sum_{I: |I||J| > 1} T_{m_I} (f_I\chi_J) \|_{L^1(CJ^c)} 
\lesssim 2^j |J|.$$
We now adapt the arguments in the previous section.  By dilation and translation invariance we may set $J = [0,1]$.  Define $\varphi$ as before, and let $r_I$ be the multipliers
$$ r_I = m_I - m_I * \varphi.$$
Then we have $T_{m_I} (f_I \chi_J) = T_{r_I} (f_I \chi_J)$ on $(CJ)^c$, and it suffices to show that
$$ \| \sum_{I: |I| > 1} T_{r_I} (f_I\chi_J) \|_{L^1(CJ^c)} 
\lesssim 2^j.$$
By H\"older as before, it suffices to show the global weighted $L^2$ estimate
$$ \| x \sum_{I: |I| > 1} T_{r_I} (f_I\chi_J)(x) \|_2
\lesssim 2^{j}.$$
By Plancherel, this becomes
\be{t} \| \sum_{I: |I| > 1} (r_I \widehat{f_I\chi_J})' \|_2
\lesssim 2^{j}.
\end{equation}
The multipliers $r_I$ can be estimated as
$$ |r_I(\xi)|, |r'_I(\xi)| \lesssim |I|^{10} (M \chi_{[\xi_I-1,\xi_I+1]})^{10}.$$
The functions $\widehat{f_I\chi_J}$ can similarly be estimated as
$$ \| \widehat{f_I\chi_J} \|_{C^1} \lesssim \|f_I \chi_J\|_1
\lesssim \| F_I * \phi_I \|_{L^1([0,1])}.$$
From the positivity of $F_I$ we have
$$ F_I * \phi_I(x) \lesssim |I|^{-10} F_I * \phi_I(0)$$
and so we thus have
$$ \| \widehat{f_I\chi_J} \|_{C^1} \lesssim |I|^{-10} (F_I * \phi_I)(0).$$
We can thus majorize the left-hand side of \eqref{t} by
$$ \| \sum_{I: |I| > 1} (F_I * \phi_I)(0) (M \chi_{[\xi_I-1,\xi_I+1]})^{10} \|_2.$$
By the Fefferman-Stein vector-valued maximal inequality \cite{feff:max}, \eqref{overlap}, and Cauchy-Schwarz as in the previous section, this is majorized by
$$ N^{1/2} (\sum_I (F_I * \phi_I)(0)^2)^{1/2} = N^{1/2} G(0) = 2^j$$
as desired.  This completes the proof of \eqref{dist} and hence \eqref{l2-est}.
\endprf

\section{Remarks on multipliers of bounded $s$-variation}\label{remarks}

Let $1 \leq s < \infty$.  For any function $f$ supported on an interval $[a,b]$, we define the $s$-variation of $f$ to be the supremum of the quantity
$$ (\sum_{i=0}^N |f(a_{i+1}) - f(a_i)|^s)^{1/s}$$
where $a = a_0 < a_1 < \ldots < a_N = b$ ranges over all partitions of $[a,b]$ of arbitrary length.  We say that a multiplier $T_m$ is a \emph{$V_s$ multiplier} if the frequency component $m_j$ have bounded $s$-variation uniformly in $j$.

Clearly the Marcinkiewicz class is the same as the $V_1$ class, but for $s > 1$ the $V_s$ class contains multipliers not covered by the Marcinkiewicz multiplier theorem.

In \cite{crdefs} it was shown that the $V_s$ class was contained in the $R_2$ class for $s < 2$.  In particular, they showed that $V_s$ multipliers were bounded on $L^p$ for $1 < p < \infty$ and $s < 2$.  From Theorem \ref{marcin-thm} and Theorem \ref{rough-thm}, we have the sharp endpoint version of this result when $s < 2$:

\begin{corollary}  Let $1 \leq s < 2$.  Then the statements of Theorem \ref{marcin-thm} (both positive and negative) continue to hold when the Marcinkiewicz class is replaced by the $V_s$ class.
\end{corollary}

Now consider the case $s > 2$.  By complex interpolation it was shown in \cite{crdefs} (see also earlier work in \cite{hirschman:multiplier}) that $V_s$ multipliers were bounded in $L^p$ when
$$|\frac{1}{2} - \frac{1}{p}| < \frac{1}{s}.$$
From the study \cite{hirschman:multiplier} of the multipliers \eqref{mab} it is known that this restriction on $p$ is sharp up to endpoints.  However, the endpoint problem remains unresolved.  The most interesting case is when $s=2$.  From the counterexamples in Section \ref{negative} we see that negative results in Theorem \ref{marcin-thm} hold for $V_2$ multipliers, and so one may conjecture that these multipliers also map both $H^1$ and $L \log^{1/2} L$ locally to $L^{1,\infty}$.  If this were true, then for $s > 2$ the $V_s$ multiplier 
class would map $L^p$ to $L^{p,p'}$ when $\frac{1}{p} = \frac{1}{s} + \frac{1}{2}$ by complex interpolation (cf. \cite{feff:thesis}).  However, we have been unable to prove these estimates using the techniques in this paper.  A natural model case would be when the frequency components $m_j$ not only have bounded $2$-variation, but have the stronger property of H\"older continuity of order $1/2$ uniformly in $j$.  (In \cite{crdefs} it was shown that a general function of bounded $2$-variation can be transformed into a H\"older continuous function of order $1/2$ by a change of variables).

In \cite{crdefs} $V_2$ multipliers were shown to be bounded on $L^p$ for all $1 < p < \infty$.  By going through their argument carefully one can show that the $L^p$ operator norm grows like $O(1/(p-1)^C)$ for some constant $C$ as $p \to 0$, so by extrapolation they map $L \log^C L$ to $L^1$ locally for some sufficiently large $C$.  However these results are far from best possible.

\section{Appendix: proof of Proposition \ref{cont}}

We now prove Proposition \ref{cont} when $f$ is in $L \log^{1/2} L([-C,C])$ and has mean zero.  

It will be convenient to move to the dyadic setting\footnote{We remark that Zygmund's original proof of \eqref{lacun} also proceeded via a dyadic model.} as we will need to perform a delicate induction shortly.  Accordingly, we introduce the Haar wavelet system
$$ \psi_I = |I|^{-1/2} (\chi_{I_l} - \chi_{I_r})$$
defined for all dyadic intervals $I$ in $[0,1]$, where $I_l$, $I_r$ are the left and right halves of $I$ respectively. 

The dyadic analogue of Proposition \ref{cont} is

\begin{proposition}\label{gen}  Let $f$ be a function on $[0,1]$ such that
$$ \int |f| \log^{1/2}(2 + |f|) \lesssim 1.$$
Then for each integer $j \geq 0$ we may find a non-negative function $f_j$ supported on $[0,1]$ such that
\be{mean}
|\langle f, \psi_I \rangle| \leq |I|^{-1/2} \int_I f_j
\end{equation}
for all $j \geq 0$ and dyadic intervals $I \subset [0,1]$ of length $2^{-j}$, and that
\be{square}
\| (\sum_{j \geq 0} |f_j|^2)^{1/2} \|_1 \lesssim 1.
\end{equation}
\end{proposition}

We now show that Proposition \ref{gen} implies Proposition \ref{cont}.  The idea is to use an averaging over translations to smooth out the dyadic singularities of the Haar wavelet system.  

Let $f$ be as in Proposition \ref{cont}; we may assume that $f$ is supported on the interval $[1/3,2/3]$.  For negative $j$, we define $F_j = |\tilde \Delta_j f|$ as in the $H^1$ theory, so that \eqref{fj-support} holds as before. From the mean zero condition of $f$ we see that $\|F_j\|_1 \lesssim 2^j$, so the contribution of these $j$ to \eqref{fj-norm} is acceptable.

For all $-1/3 \leq \theta \leq 1/3$, let $f^\theta$ denote the translated function $f^\theta(x) = f(x-\theta)$.  These functions all satisfy the requirements of Proposition \ref{gen}, with the associated functions $f^\theta_j$.  We now define $F_j$ for $j \geq 0$ by
$$ F_j(x) = \sum_{k \geq 0} 2^{-|j-k|/2} \int_{-1/3}^{1/3}  f^\theta_k(x+\theta)\ d\theta.$$

We now verify \eqref{fj-support}.  Fix $x \in [0,1]$ and $j \geq 0$.  We say that a number $-1/3 \leq \theta \leq 1/3$ is \emph{normal} with respect to $x$ and $j$ if
$$ \dist(x+\theta, 2^{-k}\Z) \geq \frac{1}{100} 2^{-|j-k|/10} 2^{-k}$$
for all integers $0 \leq k \leq j$.

Let $\Theta_{x,j}$ denote the set of all normal $\theta$; it is easy to see that $|\Theta_{x,j}| \sim 1$.  Let $\theta$ be any element of $\Theta_{x,j}$. We compute
\bas
|\Delta_j f(x)| &= |\Delta_j f^\theta(x+\theta)\\
&=|\sum_I \langle f^\theta, \psi_I \rangle \Delta_j \psi_I(x+\theta)|\\
&\leq \sum_k \sum_{I: |I| = 2^{-k}} (\int_I f^\theta_k) |I|^{-1/2} |\Delta_j \psi_I(x+\theta)|
\end{align*}
If $k \geq j$, then a computation shows that
$$ |I|^{-1/2} |\Delta_j \psi_I(x+\theta)| \lesssim 2^{2j-k} (1 + 2^k \dist(x+\theta,I))^{-100} \lesssim 2^{-|k-j|/2} 2^j 
(1 + 2^j \dist(x+\theta,I))^{-3/2}$$
and thus that
$$ \sum_{I: |I| = 2^{-k}} (\int_I f^\theta_j) |I|^{-1/2} |\Delta_j \psi_I(x+\theta)| \lesssim 2^{-|k-j|/2} f^\theta_k * \phi_j.$$
Now suppose that $k < j$.  A computation using the normality of $\theta$ shows that
$$ |I|^{-1/2} |\Delta_j \psi_I(x+\theta)| \lesssim 2^{-100|k-j|} 2^j (1 + 2^j \dist(x+\theta, I))^{-100}$$
and hence that
$$ \sum_{I: |I| = 2^{-k}} (\int_I f^\theta_j) |I|^{-1/2} |\Delta_j \psi_I(x+\theta)| \lesssim 2^{-|k-j|/2} f^\theta_k * \phi_j.$$
Combining these estimates and then averaging over $\Theta_{x,j}$ we obtain \eqref{fj-support} as desired.

Now we show \eqref{fj-norm} for the non-negative $j$.  From Young's inequality and Minkowski's inequality we see the pointwise estimate
\bas
(\sum_j |F_j(x)|^2)^{1/2} &\lesssim
(\sum_k |\int_{-1/3}^{1/3}  f^\theta_k(x+\theta)^2\ d\theta|^2)^{1/2}\\
&\leq 
\int_{-1/3}^{1/3} (\sum_k f^\theta_k(x+\theta)^2)^{1/2}\ d\theta.
\end{align*}
The claim then follows from Fubini's theorem and \eqref{square}.

It remains to prove Proposition \ref{gen}.  To do this, we first reduce to the case when $f$ is a characteristic function.  More precisely, we shall show

\begin{proposition}\label{char}  Let $N \geq 0$ be an integer, $I_0$ be a dyadic interval, and let $\I_0$ be the collection of all dyadic intervals in $I_0$ of side-length at least $2^{-N} |I_0|$.  Let $E$ be the union of some intervals in $\I$.  Then for each dyadic interval $I \subseteq I_0$ of length at least $2^{-N} |I_0|$, we may find a non-negative function $f_I$ supported on $I$ such that
\be{mean-2}
|\langle \chi_E, \psi_I \rangle| \leq |I|^{-1/2} \| f_I \|_1
\end{equation}
for all such $I$, and that\footnote{If $|E| = 0$, we adopt the convention that $|E| \log(2 + |I_0|/|E|)^{1/2} = 0$.}
\be{square-2}
\| (\sum_{I \in \I_0} |f_I|^2)^{1/2} \|_1 \leq A |E| \log(2 + |I_0|/|E|)^{1/2}
\end{equation}
for some absolute constant $A$.
\end{proposition}

Indeed, by setting $I_0 = [0,1]$ and $N \to \infty$, we see that Proposition \ref{char} immediately implies Proposition \ref{gen} for the $L \log^{1/2} L$-normalized functions $|E|^{-1} \log(1/|E|)^{-1/2} \chi_E$ for any set $E$ with measure $0 < |E| \ll 1$.  A general $L \log^{1/2} L$ function can be written as a convex linear combination of such functions (see e.g. \cite{tao:factor}), so the general case of Proposition \ref{gen} obtains (observing that the $L^1(l^2)$ space appearing in \eqref{square} is a Banach space).

It remains to prove Proposition \ref{char}.  This shall be done by induction on $N$.  Clearly the claim is true for $N=0$ simply by setting $f_{I_0} = \chi_E$.  We warn the reader in advance that the inductive nature of the argument will require some delicate estimates in which one cannot afford to lose constant factors in the main terms.

Now fix $N > 0$, $m > 0$, $I_0$, $E$, and suppose the claim holds for all smaller values of $N$.  We may rescale $I_0$ to be the unit interval $[0,1]$.

Let $0 < \eps \ll 1$ be a small absolute constant to be chosen later.
We first prove the claim in the easy case $|E| \geq \eps$.  In this case we set
$$ f_I = |I|^{-1/2} |\langle \chi_E, \psi_I \rangle| \chi_I.$$
The estimate \eqref{mean-2} is trivial.  To verify \eqref{square-2}, we use H\"older's inequality and the orthonormal nature of the Haar basis:
\bas
\| (\sum_{I \in \I_0} |f_I|^2)^{1/2} \|_1 &\leq
\| (\sum_{I \in \I_0} |f_I|^2)^{1/2} \|_2\\
&= (\sum_{I \in \I_0} |\langle \chi_E, \psi_I \rangle|^2)^{1/2}\\
&\leq \| \chi_E \|_2\\
&\lesssim |E| \log(2 + 1/|E|)^{1/2}
\end{align*}
as desired (if $A$ is sufficiently large depending on $\eps$.

Now suppose $|E| < \eps$.  Let $\I$ denote the set of all intervals $I \in \I_0$ such that
\be{triad}
\eps |E||I| \leq |E \cap I| \geq 2|E| |I|.
\end{equation}
holds, where $0 < \eps \ll 1$ is an absolute constant to be chosen later.  Let $\J$ denote the set of all intervals not in $\I$ which are maximal with respect to set inclusion.  From our assumptions on $E$ we see that $\J$ is a partition of $[0,1]$ into disjoint intervals, and each interval $J \in \J$ satisfies
$$ 2^{-N} < |J| < 1.$$

Let $J$ be any element of $\J$.  From the induction hypothesis we can associate a function $f_I$
to each $I \in \I_0$, $I \subseteq J$ such that
$$ \langle \chi_E, \psi_I \rangle = \langle \chi_{E \cap J}, \psi_I \rangle \leq |I|^{-1/2} \int_I f_I$$
for all such $I$, and 
\be{low}
\| F_J \|_1 \leq A |E \cap J| \log(2 + |J|/|E \cap J|)^{1/2},
\end{equation}
where we have written $F_J$ for the function
$$ F_J = (\sum_{I \in \I_0: I \subseteq J} |f_I|^2)^{1/2}.$$

We have now defined the $f_I$ for all intervals contained in one of the intervals $J \in \J$.  It remains to assign functions $f_I$ to the intervals $I$ in $\I$.  

Let $\I^*$ denote those intervals $I$ in $\I$ such that $|E \cap I| > 0$.  We will set $f_I = 0$ for all $I \in \I \backslash \I^*$; note that \eqref{mean-2} holds vacuously for these $I$.  For $I \in \I^*$, we define $f_I$ by the formula
$$ f_I = |I|^{1/2} |\langle \chi_E, \psi_I \rangle| \sum_{J \in \J: J \subset I} \frac{|E \cap J|}{|E \cap I|} \frac{F_J}{\|F_J\|_1}.$$
Since $I$ is the union of the intervals $J \in \J$ contained inside it, we see that 
$$ \|f_I\|_1 = |I|^{1/2} |\langle \chi_E, \psi_I \rangle| \sum_{J \in \J: J \subset I} \frac{|E \cap J|}{|E \cap I|} = |I|^{1/2} |\langle \chi_E, \psi_I \rangle|$$
so that \eqref{mean-2} holds for these $I$.  

We now verify \eqref{square-2}.  For any $J \in \J$ and $x \in J$, we have
\bas
\sum_{I \in \I_0} |f_I(x)|^2
&= \sum_{I \in \I_0: I \subseteq J} |f_I(x)|^2 + 
\sum_{I \in \I^*: I \supset J} |f_I(x)|^2)^{1/2} \\
&= F_J(x)^2 + 
\sum_{I \in \I^*: I \supset J } |I| |\langle \chi_E, \psi_I \rangle|^2 \frac{|E \cap J|^2}{|E \cap I|^2} \frac{F^2_J(x)}{\|F_J\|^2_1}\\
&=
\frac{F_J(x)^2}{\|F_J\|^2_1} ( \|F_J\|_1^2 + \sum_{I \in \I^*: I \supset J} |I| \frac{|E \cap J|^2}{|E \cap I|^2} |\langle \chi_E, \psi_I \rangle|^2 ).
\end{align*}
Taking the square root of this and integrating, we obtain
\be{id}
\| (\sum_{I \in \I_0} |f_I|^2)^{1/2} \|_1
= 
\sum_{J \in \J} 
(\|F_J\|_1^2  + \sum_{I \in \I^*: I \supset J} |I| \frac{|E \cap J|^2}{|E \cap I|^2} |\langle \chi_E, \psi_I \rangle|^2)^{1/2}.
\end{equation}
Now define the function 
$$ g = \sum_{J \in \J} |E \cap J| \frac{\chi_J}{|J|}.$$
For all $I \in \I^*$ we see that $\psi_I$ is constant on intervals in $\J$, and hence that $\langle g, \psi_I \rangle = \langle \chi_E, \psi_I\rangle$.  Thus
\be{id2}
\eqref{id} = 
\sum_{J \in \J} 
(\|F_J\|_1^2  + \sum_{I \in \I^*: I \supset J} |I| \frac{|E \cap J|^2}{|E \cap I|^2} |\langle g, \psi_I \rangle|^2)^{1/2}.
\end{equation}
For future reference we observe from the construction of $\J$ and $g$ that $\|g\|_1 = |E|$ and $\|g\|_\infty \leq 4 |E|$, hence
\be{lump} \sum_{I \in \I^*} |\langle g, \psi_I \rangle|^2 \leq \|g\|_2^2 \leq \|g\|_1 \|g\|_\infty \lesssim |E|^2.
\end{equation}

To estimate \eqref{id2}, we define 
\bas
\J_1 &= \{ J \in \J: 2|E| |J| \leq |E \cap J| \leq 4|E| |J| \} \\
\J_2 &= \{ J \in \J: |E|^{10} |J| \leq |E \cap J| \leq \eps|E| |J| \} \\
\J_3 &= \{ J \in \J: |E \cap J| < |J| |E|^{10} \};
\end{align*}
note from \eqref{triad} and the construction of $\J$ that $\J = \J_1 \cup \J_2 \cup \J_3$.  Thus \eqref{id} is the sum of 
\be{j1}
\sum_{J \in \J_1 \cup \J_2} 
(\|F_J\|_1^2  + \sum_{I \in \I^*: I \supset J} |I| \frac{|E \cap J|^2}{|E \cap I|^2} |\langle g, \psi_I \rangle|^2)^{1/2}.
\end{equation}
and
\be{j3}
\sum_{J \in \J_3} 
(\|F_J\|_1^2  + \sum_{I \in \I^*: I \supset J} |I| \frac{|E \cap J|^2}{|E \cap I|^2} |\langle g, \psi_I \rangle|^2)^{1/2}.
\end{equation}

We first consider \eqref{j3}, the contribution of the very sparsely occupied intervals.  In this case we use crude estimates.  From the estimate $(a^2 + b)^{1/2} \leq a+b^{1/2}$ we have
$$
\eqref{j3} \leq \sum_{J \in \J_3} \|F_J\|_1 + \sum_{J \in \J_3}
(\sum_{I \in \I^*: I \supset J} |I| \frac{|E \cap J|^2}{|E \cap I|^2} 
|\langle g, \psi_I \rangle|^2)^{1/2}
$$
To estimate the first term, we observe from \eqref{low} that
$$
\| F_J \|_1 \lesssim A |E|^{10} |J| \log(1/|E|)^{1/2}$$
and so
$$
\sum_{J \in \J_3} \|F_J\|_1 \lesssim A |E|^{10} \log(1/|E|)^{1/2} \lesssim A |E|^9$$
since we of course have
\be{j3-est}
\sum_{J \in \J_3} |J| \leq 1
\end{equation}
To estimate the second term, we use Cauchy-Schwarz and \eqref{j3-est}, to obtain
$$
\eqref{j3} \leq C A |E|^9 + (\sum_{J \in \J_3}
|J|^{-1} \sum_{I \in \I^*: I \supset J} |I| \frac{|E \cap J|^2}{|E \cap I|^2} 
|\langle g, \psi_I \rangle|^2)^{1/2}.
$$
Using the estimate $|J|^{-1} |E \cap J| \leq |E|^{10}$, and then interchanging summations, we obtain
$$
\eqref{j3} \leq C A |E|^9 + (\sum_{I \in \I^*} \sum_{J \in \J: J \subset I}
|E|^{10} |I| \frac{|E \cap J|}{|E \cap I|^2} 
|\langle g, \psi_I \rangle|^2)^{1/2}.
$$
Performing the $J$ summation, this becomes
$$
\eqref{j3} \leq C A |E|^9 + |E|^5 (\sum_{I \in \I^*} \frac{|I|}{|E \cap I|} 
|\langle g, \psi_I \rangle|^2)^{1/2}.
$$
Applying \eqref{triad} and then \eqref{lump} we thus obtain
\be{j3a-est}
\eqref{j3} \leq C A |E|^9 + |E|^5 ( |E|^{-1} |E|^2)^{1/2} \lesssim A |E|^2.
\end{equation}

Now we turn to the more interesting term \eqref{j1}.  From \eqref{low} we have
$$ \eqref{j1} \leq
\sum_{J \in \J_1 \cup \J_2} 
((A |E \cap J| \log(2 + |J|/|E \cap J|)^{1/2})^2  + \sum_{I \in \I^*: I \supset J} |I| \frac{|E \cap J|^2}{|E \cap I|^2} |\langle g, \psi_I \rangle|^2)^{1/2}.$$
Using the inequality 
$$ \sqrt{a^2 + b} \leq \sqrt{a^2 + b + \frac{b^2}{4a^2}} = a + \frac{b}{2a},$$
for $a, b > 0$, we thus have
$$
\eqref{j1} \leq \eqref{j4} + \eqref{j5}$$
where \eqref{j4} and \eqref{j5} are given by
\be{j4}
\sum_{J \in \J_1 \cup \J_2} A |E \cap J| \log(2 + |J|/|E \cap J|)^{1/2}
\end{equation}
and
\be{j5}
\sum_{J \in \J_1 \cup \J_2} \frac{1}{2A |E \cap J| \log(2 + |J|/|E \cap J|)^{1/2}}
\sum_{I \in \I^*: I \supset J} |I| \frac{|E \cap J|^2}{|E \cap I|^2} |\langle g, \psi_I \rangle|^2.
\end{equation}

Let us first estimate the error term \eqref{j5}.  Since $J \in \J_1 \cup \J_2$, we see that
$$
\log(2 + |J|/|E \cap J|)^{1/2} \sim \log(1/|E|)^{1/2}.$$
Applying this, re-arranging the summation, and simplifying, we obtain
$$
\eqref{j5} \lesssim
\log(1/|E|)^{-1/2} \sum_{I \in \I^*} \sum_{J \in \J: J \subset I} 
|I| \frac{|E \cap J|}{|E \cap I|^2} |\langle g, \psi_I \rangle|^2.$$
Performing the $J$ summation, we obtain
$$
\eqref{j5} \lesssim
\log(1/|E|)^{-1/2} \sum_{I \in \I^*} \frac{|I|}{|E \cap I|} |\langle g, \psi_I \rangle|^2.$$
From \eqref{triad} and \eqref{lump} we thus have
\be{j5-est}
\eqref{j5} \lesssim |E| \log(1/|E|)^{-1/2}.
\end{equation}

It remains to treat \eqref{j4}, which is the main term.  We split this as
$\eqref{j4} = \eqref{j6} - \eqref{j7} + \eqref{j8}$, where \eqref{j6}, \eqref{j7}, \eqref{j8} are given by
\be{j6}
\sum_{J \in \J_1 \cup \J_2} A |E \cap J| \log(2 + 1/|E|)^{1/2}
\end{equation}
\be{j7}
\sum_{J \in \J_1} A |E \cap J| (\log(2 + 1/|E|)^{1/2} - \log(2 + |J|/|E \cap J|)^{1/2})
\end{equation}
\be{j8}
\sum_{J \in \J_2} A |E \cap J| (\log(2 + |J|/|E \cap J|)^{1/2} - \log(2 + 1/|E|)^{1/2}).
\end{equation}
Note that \eqref{j6}, \eqref{j7}, \eqref{j8} are all non-negative.  We can estimate \eqref{j6} by
$$ \eqref{j6} \leq A |E| \log(2 + 1/|E|)^{1/2}$$
which is exactly the quantity needed for the induction hypothesis.  Collecting all the terms and using \eqref{j3a-est}, \eqref{j5-est}, we see that we have to show that
\be{tight}
 \eqref{j7} \geq \eqref{j8} + C A |E|^2 + C |E| \log(1/|E|)^{-1/2}.
\end{equation}
We thus seek good lower bounds on \eqref{j7} and good upper bounds on \eqref{j8}.

We first deal with \eqref{j7}.  We may write this as
$$\eqref{j7} = A\sum_{J \in \J_1} |E \cap J| 
\frac{\log(2 + 1/|E|) - \log(2 + |J|/|E \cap J|)}
{(\log(2 + 1/|E|)^{1/2} + \log(2 + |J|/|E \cap J|)^{1/2}}.$$
Both terms in the denominator are comparable to $\log(1/|E|)^{1/2}$, while the numerator is bounded from below by
$$ \log(2 + 1/|E|) - \log(2 + 1/{2|E|}) \sim 1.$$
Thus we have
$$ \eqref{j7} \sim A \log(1/|E|)^{1/2} \sum_{J \in \J_1} |E \cap J|.$$
To obtain lower bounds for this, we observe that
$$ \sum_{J \in \J_1} |E \cap J| = |E| - \sum_{J \in \J_2 \cup \J_3} |E \cap J|$$
and
$$ \sum_{J \in \J_2 \cup \J_3} |E \cap J| \leq \sum_{J \in \J} \eps |E| |J| = \eps |E|.$$
Thus
$$ \eqref{j7} \gtrsim A |E| \log(1/|E|)^{-1/2}.$$
Now we attend to \eqref{j8}.  As before, we may write
$$ \eqref{j8} =
A\sum_{J \in \J_1} |E \cap J| 
\frac{\log(2 + |J|/|E \cap J|) - \log(2 + 1/|E|)}
{(\log(2 + 1/|E|)^{1/2} + \log(2 + |J|/|E \cap J|)^{1/2}}.$$
Again, the denominator is comparable to $\log(1/|E|)^{1/2}$, while the numerator is comparable to $\log(|E||J|/|E \cap J|)$.  Thus
$$ \eqref{j8} \lesssim A \log(1/|E|)^{-1/2}
\sum_{J \in \J: |E \cap J| \leq \eps |E| |J|} |E \cap J| \log(|E||J|/|E \cap J|).$$
We estimate this dyadically as
\bas
\eqref{j8} &\lesssim A \log(1/|E|)^{-1/2}
\sum_{k: 2^{-k} \lesssim \eps} \sum_{J \in \J: |E \cap J| \sim 2^{-k} |E| |J|} |E \cap J| \log(|E||J|/|E \cap J|)\\
&\lesssim A \log(1/|E|)^{-1/2} \sum_{k: 2^{-k} \lesssim \eps}
\sum_{J \in \J} 2^{-k} |E| |J| k\\
&\lesssim A |E| \log(1/|E|)^{-1/2} \sum_{k: 2^{-k} \lesssim \eps} 2^{-k} k\\
&\lesssim A |E| \log(1/|E|)^{-1/2} \sum_{k: 2^{-k} \lesssim \eps} 2^{-k/2}\\
&\lesssim A \eps^{1/2} |E| \log(1/|E|)^{-1/2}
\end{align*}

Thus \eqref{tight} resolves to
$$ C^{-1} A |E| \log(1/|E|)^{-1/2} \geq
C A \eps^{1/2} |E| \log(1/|E|)^{-1/2} + C A |E|^2 + C |E| \log(1/|E|)^{-1/2},$$
and this is achieved if $\eps$ is chosen sufficiently small (recall that $|E| \leq \eps$), and then $A$ is chosen sufficiently large depending on $\eps$.

\endprf


\begin{thebibliography}{10}

\bibitem{bergh:interp}
J. Bergh, J. L\"ofstr\"om, \emph{Interpolation Spaces: An Introduction},
Springer-Verlag, 1976.

\bibitem{crdefs}
R. R. Coifman, J. L. Rubio de Francia, S. Semmes, \emph{Multiplicateurs de Fourier dans $L^p(R)$ et estimations quadratiques}, C. R. Acad Sci. Paris \textbf{306} (1988), 351-354.

\bibitem{feff:thesis}
C. Fefferman, \emph{Inequalities for strongly
singular convolution operators}, Acta Math. \textbf{ 124}
(1970), 9--36.
 
\bibitem{feff:max}
C. Fefferman and E.~M. Stein, \emph{Some maximal inequalities}, Amer.
J. Math. \textbf{93} (1971): 107--115.

\bibitem{hirschman:multiplier}
I.I. Hirschman, \emph{On multiplier transformations}, Duke Math. J. 26 (1959), 221-254.

\bibitem{oneil}
R. O'Neil, \emph{Convolution operators and $L(p,q)$ spaces}, Duke Math. J. \textbf{30}, 129-143 (1963).

\bibitem{sagher}
Y. Sagher, \emph{On analytic families of operators}, Israel J. Math.
\textbf{7} (1969), 350--356.

\bibitem{seeger:compact}
A. Seeger, \emph{Endpoint estimates for multiplier transformations on
compact manifolds}, Indiana Math. J. \textbf{40} (1991): 471--533.

\bibitem{ts:l12}
A. Seeger, T. Tao, \emph{Sharp Lorentz space estimates for rough operators}, submitted, Math. Annalen.  {\tt math.CA/9912098}

\bibitem{stein:small}
E.~M. Stein, \emph{Singular Integrals and Differentiability
Properties of Functions}, Princeton University Press, 1970.

\bibitem{stein:large}
E.~M. Stein, \emph{Harmonic Analysis}, Princeton University Press, 1993.

\bibitem{tao:factor} T. Tao, \emph{A converse extrapolation theorem for translation invariant operators}, submitted, J. Funct. Anal.  {\tt math.FA/9912001}

\bibitem{yano} S. Yano, \emph{Notes on Fourier analysis. XXIX. An extrapolation theorem. }
J. Math. Soc. Japan \textbf{3}, (1951). 296--305.

\bibitem{zygmund}  A. Zygmund, \emph{Trigonometric series. Vol. I, II.}, Reprint of the 1979 edition. Cambridge Mathematical Library. Cambridge University
Press, Cambridge-New York, 1988.

\end{thebibliography}
\end{document}